\documentstyle{article}
\newtheorem{dfn}{Definition}
\newtheorem{prop}[dfn]{Proposition}
\newtheorem{thm}[dfn]{Theorem}

\begin{document}
\title{The multiple gamma function and its $q$-analogue
    \thanks{To appear in the proceedings of the workshop
      `` Quantum groups and Quantum spaces'' at the Banach center}}  
\author{Kimio UENO and Michitomo NISHIZAWA\\
 Department of Mathematics\\
 School of Science and Engineering\\
 Waseda University}
\maketitle

\begin{abstract}
We give an asymptotic expansion ({\it the higher Stirling formula})
and an infinite product representation ({\it the Weierstrass product
formula}) of the Vign\'{e}ras multiple gamma function by considering
the classical limit of the multiple $q$-gamma function.
\end{abstract}


\section{Introduction}
The multiple gamma function was introduced by Barnes. It is defined
to be an infinite product regularized by the multiple Hurwitz zeta-functions 
\cite{bar1}, \cite{bar2}, \cite{bar3}, \cite{bar4}. After his discovery, 
many mathematicians have studied this function:
Hardy \cite{har1}, \cite{har2} studied
this function from his viewpoint of the theory of elliptic functions,
and Shintani \cite{shi1}, \cite{shi2} applied it to the study on 
the Kronecker limit formula for zeta-functions attached to algebraic 
fields.\par
 In the end of 70's, Vign\'{e}ras \cite{vig} redefined a multiple
gamma function to be a function satisfying the generalized Bohr-Morellup
theorem, Furthermore, Vign\'{e}ras \cite{vig}, Voros \cite{vor} and Kurokawa 
\cite{kur1}, \cite{kur2}, \cite{kur3}, \cite{kur4} showed that it plays 
an essential role to express the Selberg zeta-function and the determinant 
of Laplacians.\par
As we can see from these studies, the multiple gamma function is a 
fundamental function for the analytic number theory: See also \cite{kur5},
\cite{man}. However we do not think that the theory of the multiple gamma
functions has been fully explored.\par
 On the other hand, the second author of this article introduced a
$q$-analogue of the Vign\'{e}ras multiple gamma functions and showed it to be 
characterized by a $q$-analogue of the generalized Bohr-Morellup 
theorem \cite{vig}.\par
 In this article, we will establish an asymptotic expansion formula 
({\it the higher Stirling formula}) and an infinite product representation 
({\it the Weierstrass product formula}) of the Vign\'{e}ras multiple gamma 
function by considering the classical limit of the multiple 
$q$-gamma functions. In order to get these results, we will use the method 
developed in \cite{un1}. Namely, by making use of the Euler-MacLaurin summation
formula, we derive the Euler-MacLaurin expansion of the multiple $q$-gamma 
function. Taking the classical limit, we lead to the 
Euler-MacLaurin expansion of the Vig\'{n}eras multiple gamma function.
The higher Stirling formula and the Weierstrass product formula
follow from this expansion formula \par
The details of the proof will be published in the forthcoming paper
\cite{un2}.


\section{A survey of the multiple gamma function and the multiple 
$q$-gamma function}
\subsection{The gamma function}
The following are well-known facts in the classical analysis: 
The Bohr-Morellup theorem says that the gamma function $\Gamma(z)$ 
is characterized by the three conditions, 
   \begin{enumerate}
      \item $\Gamma(z+1)=z \Gamma(z)$,
      \item $\Gamma(1)=1$,
      \item $\frac{d^2}{dz^2}\log\Gamma(z+1)\geq0
        \mbox{ for } z\geq0.$
   \end{enumerate}
The gamma function is meromorphic on {\bf C}, and has an infinite product
representation 
   \begin{equation}
     \Gamma(z+1)=e^{-\gamma x}
       \prod_{n=1}^{\infty}
        \left\{
        \left(1+\frac{z}{n}\right)^{-1}
        e^{\frac{z}{n}}
        \right\},
   \end{equation}
where $\gamma$ is the Euler Constant. This is usually called the Weierstrass
product formula.\par
The gamma function has an asymptotic expansion, which is called the 
Stirling formula,
  \begin{eqnarray*}
    & &\log \Gamma(z+1) \sim
      \left(z+\frac{1}{2}\right)\log(z+1) - (z+1) - \zeta'(0)\\
    & &\qquad + \sum_{r=1}^{\infty} \frac{B_{2r}}{[2r]_{2}}
    \frac{1}{(z+1)^{2r-1}},
  \end{eqnarray*}
as $z\to\infty$ in the sector
$\Delta_{\delta}:= \{z\in{\bf C}||\arg z|<\pi-\delta\}$  $(0<\delta<\pi),$
where 
     $$\frac{z e^{tz}}{e^{z}-1}
        = \sum_{n=0}^{\infty} \frac{B_{n}(t)}{n!} z^{n}$$
$B_{k}=B_{k}(0)\quad(\mbox{the Bernoulli number})$,  $\zeta(s)$ 
is the Riemann zeta-function, $\zeta'(s)= \frac{d}{ds}\zeta(s)$ and 
$[x]_{r}=x(x-1)\cdots (x-r+1)$. Note that $\zeta'(0)=-\log \sqrt{2\pi}$.

\subsection{The Barnes $G$-function}
Barnes \cite{bar1} introduced  the function $G(z)$ which satisfies
   \begin{enumerate}
      \item $G(z+1)= \Gamma(z) G(z)$,
      \item $G(1)=1$,
      \item $\frac{d^3}{dz^3}\log G(z+1)\geq0
        \mbox{ for } z\geq0,$
   \end{enumerate}
and he called this ``$G$-function''. He proved that the $G$-function has an 
infinite product representation.
  $$ G(z+1)
      = e^{-z\zeta'(0)-\frac{z^{2}}{2}\gamma
      - \frac{z^{2}+z}{2}}
      \prod_{k=1}^{\infty} \left\{\left(
        1+\frac{z}{k}\right)^{k}
        \exp\left(-z+\frac{z^{2}}{2k}\right)\right\}.$$
and an asymptotic expansion 
  \begin{equation} 
    \log G(z+1)\sim\left(\frac{z^{2}}{2}-\frac{1}{12}\right)\log(z+1)
      -\frac{3}{4}z^{2}-\frac{z}{2}+\frac{1}{3}+z\zeta'(0)
      -\log A + O(\frac{1}{z})
  \label{eqn:asG}\end{equation}
as $z\to\infty$ in the sector $\Delta_{\delta}$ where $A$ is called the 
Kinkelin constant. Voros showed this constant can be written with the first 
derivative of the Riemann zeta-function
(cf \cite{vor}, \cite{var})
  $$ \log A = -\zeta'(-1)+\frac{1}{12}.$$

\subsection{The Vign\'{e}ras multiple gamma function} 
As a generalization of the gamma function and the $G$-function, Vign\'{e}ras 
\cite{vig} introduced a hierarchy of functions which satisfy
  \begin{eqnarray*}
      & & 1. \quad G_{n}(z+1)=G_{n-1}(z) G_{n}(z), \\
      & & 2. \quad G_{n}(1)=1, \\
      & & 3. \quad \frac{d^{n+1}}{dz^{n+1}}\log G_{n}(z+1)\geq0
         \quad \mbox{ for } \quad z\geq0,\\
      & & 4. \quad G_{0}(z) =z.
  \end{eqnarray*}
and she called these functions ``the multiple gamma functions''. 
Applying Dufresnoy and Pisot's results \cite{duf}, she showed that 
these functions satisfying the above properties are unique and that 
it has an infinite product representation: 
  \begin{eqnarray}
     & &G_{n}(z+1)=\exp \left[-zE_{n}(1)+\sum_{h=1}^{n-1}
       \frac{p_{h}(z)}{h!}\left(\psi_{n-1}^{(h)}(0)-E_{n}^{(h)}(1)\right)
       \right]\nonumber\\
     & &\nonumber\\  
     & & \times \prod_{{\bf m}\in {\bf N}^{n-1} \times {\bf N}^{*}}
       \left[\left(1+\frac{z}{s({\bf m})}\right)^{(-1)^{n}}
       \exp \left\{ \sum_{l=0}^{n-1}\frac{(-1)^{n-l}}{n-l}
       \left(\frac{z}{s({\bf m})}\right)^{n-l}\right\}\right]\\
     & &\nonumber
  \label{eqn:vig}\end{eqnarray}           
where 
  \begin{eqnarray*}
    & & E_{n}(z):= \sum_{{\bf m}\in {\bf N}^{n-1} \times {\bf N}^{*}}
      \left[\left\{\sum_{l=0}^{n-1}\frac{(-1)^{n-l}}{n-l}
        \left(\frac{z}{s({\bf m})}\right)^{n-l}\right\}
        +(-1)^{n} \log \left(1+\frac{z}{s({\bf m})}\right)\right],\\
    & &\\    
    & & \psi_{n-1}^{(h)}(z):=\frac{d^{h}}{dz^{h}}\log G_{n-1}(z+1),\\    
    & &\\
    & & p_{h}(z):=1^{h}+2^{h}+\cdots+(z-1)^{h},\\
    & &\\
    & & s({\bf m}):=m_{1}+m_{2}+\cdots+m_{n}\quad \mbox{for}
        \quad {\bf m}=(m_{1},m_{2},\cdots m_{n}),
  \end{eqnarray*}
and ${\bf N}^{*}={\bf N}-\{0\}$.

\subsection{The $q$-gamma function}
Throughout this article, we suppose $0<q<1$. 
A $q$-analogue of the gamma function was defined by Jackson \cite{jac1}, 
\cite{jac2}: 
  $$\Gamma (z+1;q)=(1-q)^{-z}
        \prod_{k=1}^{\infty}\left(
            \frac{1-q^{z+k}}{1-q^{k}}
            \right)^{-1}.$$
Askey \cite{ask} pointed out that this function satisfies 
a $q$-analogues of the Bohr-Morellup theorem. Namely, $\Gamma(z;q)$ 
satisfies
    \begin{enumerate}
        \item $\Gamma (z+1;q)=[z] \Gamma(z;q),$
        \item $\Gamma (1;q)=1,$
        \item $\frac{d^{2}}{dz^{2}}\log\Gamma(z+1;q)\geq0
            \mbox{ for } z\geq0,$
    \end{enumerate}
where $[z]:=(1-q^{z})/(1-q)$. \par
As $q$ tends to $1-0$, $\Gamma(z;q)$ converges $\Gamma(z)$ uniformly with
respect to $z$. A rigorous proof of this fact  was given by Koornwinder 
\cite{koo}.\par
 Inspired by Moak's works \cite{moa}, the authors \cite{un1} derived a
representation of the $q$-gamma function    
  \begin{eqnarray}
    \log \Gamma(z:q) & = & \left(z-\frac12\right) \log \left(
        \frac{1-q^{z}}{1-q} \right)
        + \log q \int_{1}^{z} \xi \frac{q^{\xi}}{1-q^{\xi}}d\xi \nonumber\\
        & &\nonumber\\
        & + &C_{1}(q) + \frac{1}{12} \log q
        + \sum_{k=1}^{\infty} \frac{B_{2k}}{(2k)!} \left(
        \frac{\log q}{q^{z}-1} \right)^{2k-1} P_{2k-1}(q^{z})\nonumber\\
        & &\nonumber\\
        & + &R_{2m}(z;q) \label{eqn:2}
  \end{eqnarray}
where
  \begin{eqnarray*}
    & &C_{1}(q) = - \frac{1}{12}\log q -\frac{1}{12}\,\frac{\log q}{q-1}
           + \int _{0}^{\infty} \frac{\overline{B}_{2}(t)}{2}
           \left( \frac{\log q}{q^{t+1}-1} \right)^{2} q^{t+1}dt,\\
    & &\\       
    & & R_{2m}(z:q) = \int_{0}^{\infty} \frac{\overline{B}_{2m}(t)}{(2m)!}
                \left(
                \frac{\log q}{q^{t+z}-1}
                \right)^{2m} P_{2m}(q^{t+z})dz,\\
    & & \nonumber            
  \end{eqnarray*} 
and the polynomial $P_{n}(x)$ is defined by the recurrence formula
  $$P_{1}(z)=1, \quad 
    (x^{2}-x)\frac{d}{dx}P_{n}(x)+nxP_{n}(x)=P_{n+1}(x).$$ 
Each term of the formula (\ref{eqn:2}) converges uniformly as $q\to 1-0$. 
So we get another proof of the uniformity of the classical limit of 
$\log \Gamma(z;q).$

\subsection{The multiple $q$-gamma function}
Recently, one of the authors \cite{nis} construct the function 
$G_{n}(z;q)$ which satisfies a $q$-analogue of the generalized Bohr-Morellup
theorem  
    \begin{enumerate}
        \item $G_{n}(z+1;q)=G_{n-1}(z;q)G_{n}(z;q),$
        \item $G_{n}(1;q)=1,$
        \item $\frac{d^{n+1}}{dz^{n+1}}\log G_{n+1}(z+1;q)\geq0
            \mbox{ for } z\geq0,$
        \item $G_{0}(z;q)=[z].$
    \end{enumerate}
We call it ``the multiple $q$-gamma function''. It is given by the following
infinite product representation:
  \begin{eqnarray}
    & &G_{0}(z+1;q):= [z+1]\nonumber\\
    & &\nonumber\\
    & &G_{n}(z+1;q):= (1-q)^{- {z \choose n}} \prod_{k=1}^{\infty}
        \left\{
        \left(
        \frac{1-q^{z+k}}{1-q^{k}}
        \right)^{-k \choose n-1}
        (1-q^{k})^{g_{n}(z,k)}
        \right\} \label{eqn:Gn}\\
    & &\qquad \qquad \quad (n \geq 1),\nonumber
  \end{eqnarray}
where
  $$g_{n}(z;u)={z-u \choose n-1}-{-u \choose n-1}.$$ \par
In the next section, we derive a representation of the multiple $q$-gamma
function like (\ref{eqn:2}) and consider its classical limit. This limit
formula gives some important properties of the multiple gamma function.


\section{The Euler-MacLaurin expansion}
Our main tool is the Euler-MacLaurin summation formula
  \begin{eqnarray*}
    \sum_{r=M}^{N-1} f(r)  =  \int_{M}^{N}f(t)dt
        & + &  \sum_{k=1}^{n} \frac{B_{k}}{k!}
        \left\{
        f^{(k-1)}(N)-f^{(k-1)}(M)
        \right\} \nonumber\\
        & + & (-1)^{n-1} \int_{M}^{N} \frac{\overline{B}_{n}(t)}{n!}
        f^{(n)}(t)dt.
  \end{eqnarray*}
where $\overline{B}_{n}(t)=B_{n}(t-[t])$ and $[t]$ denotes an integral 
part of $t$.
We apply this to $\log G_{n}(z+1;q)$.\par
 From the infinite product representation (\ref{eqn:Gn}),  we have 
  \begin{eqnarray}
    \log G_{n}(z+1;q)
      &=& -{z \choose n}\log(1-q) - \sum_{k=1}^{\infty}
        {-k \choose n-1}\log (1-q^{z+k})\nonumber\\
      & &\nonumber\\  
      &+& \sum_{j=0}^{n-1} G_{n,j}(z)\{\sum_{k=1}^{\infty}
        k^{j}\log(1-q^{k})\}\label{eqn:31},
  \end{eqnarray}
where
  $${z-u \choose n-1} = \sum_{j=0}^{n-1}G_{n,j}(z)u^{j}.$$
By the Euler-MacLaurin summation formula, we obtain
  \begin{eqnarray}
    & & \sum_{k=1}^{\infty}{-k \choose n-1}\log(1-q^{z+k})\nonumber\\
    & &\nonumber\\
    & & \quad = \sum_{j=0}^{n-1}
      \frac{(-1)^{j}{}_{n-1}S_j}{(n-1)!}
      \sum_{r=0}^{j} \frac{(-1)^{r}j!}{(j-r)!}L_{r+2}(z+1)\nonumber\\
    & &\nonumber\\  
    & & \quad    + \sum_{r=1}^{m} \frac{B_{r}}{r!}
      \left\{ \left(
        \frac{d}{dt} \right)^{r-1} {-t \choose n-1}
      \right\}\Bigg\vert_{t=1}L_{1}\nonumber\\
    & &\nonumber\\  
    & & \quad - \sum_{r=1}^{m} \frac{B_{r}}{r!}F_{n, r-1}(z;q)
       + R_{n,m}(z;q),\label{eqn:33}
  \end{eqnarray}
where  
 \begin{eqnarray*}
    & &L_{r}(z):= \frac{Li_{r}(q^{z})}{\log^{r-1}q}\mbox{ },
      \quad L_{1}(z):= -\log(1-q^{z}), \quad
      L_{r}:=L_{r}(1),\\ 
    & &\\  
    & & Li_{r}(z):=\sum_{k=1}^{\infty}\frac{z^{k}}{k^{r}}
      \qquad (\mbox{Euler's polylogarithm}),\\
    & &\\
    & & F_{n,r-1}(z;q) := \left[\frac{d^{r-1}}{dt^{r-1}}
      \left\{ {-t \choose n-1}\log \left( \frac{1-q^{z+t}}{1-q^{z+1}}
      \right)
      \right\}\right]_{t=1},\\
    & &\\
    & &  R_{n,m}(z;q):=  \frac{(-1)^{m-1}}{m!}
         \int_{1}^{\infty}
         \left[ \overline{B}_{m}(t)
         \left\{\frac{d^{m}}{dt^{m}}
            \left\{{-t \choose n-1} \log
            \left(\frac{1-q^{z+t}}{1-q^{z+1}}\right)
            \right\}
         \right\} \right]dt,\\
   & & 
\end{eqnarray*}
and ${}_nS_j$ is the Stirling number of the first kind
defined by
 $$ \sum_{j=0}^{n}{}_nS_j u^{j}
    = u(u-1)(u-2)\cdots(u-n+1).$$
Similarly, we have
  \begin{equation}  
      \sum_{k=1}^{\infty} k^{j}\log(1-q^{k})
      = \sum_{r=0}^{j}
      \frac{(-1)^{r}j!}{(j-r)!} L_{r+2}
      + \sum_{r=1}^{j}\frac{B_{r}}{r!}
        \frac{j!}{(j+1-r)!}L_{1} + C_{j}(q), \label{eqn:38}
  \end{equation}
where  
   \begin{eqnarray*} 
     & & C_{j}(q):= -\sum_{r=1}^{n+1}\frac{B_{r}}{r!}f_{j+1,r-1}(q)\\
     & &\\
     & & \qquad \quad
         + \frac{(-1)^{n}}{(n+1)!}\int_{1}^{\infty}
         \left[ \overline{B}_{n+1}(t)
         \left\{\frac{d^{n+1}}{dt^{n+1}}
            \left\{t^{j}\log
              \left(\frac{1-q^{t}}{1-q} \right)
            \right\}
         \right\} \right]dt,\\
    & &\\
    & & f_{j+1,r-1}(q):=
         \left[\frac{d^{r-1}}{dt^{r-1}}
            \left\{t^{j} \log
            \left( \frac{1-q^{t}}{1-q} \right)
            \right\} \right]_{t=1}.\\
   & &
  \end{eqnarray*}          
Substituting  (\ref{eqn:33}) and (\ref{eqn:38}) to (\ref{eqn:31}),
we have 
  \begin{eqnarray}
    & &\log G_{n}(z+1;q)\nonumber\\
    & &\quad = ( \mbox{the terms containing } L_{r}(z) \mbox{ and }
      L_{r})\nonumber\\
    & & \qquad \quad + ( \mbox{the term which converges as } q\to1).
  \label{eqn:3*}\end{eqnarray}    
$L_{r}(z)$ and $L_{r}$ give rise to a divergent part in the expression
(\ref{eqn:3*}) as $q\to1-0$.
But, thanks to the identity
  \begin{eqnarray*}
    & & L_{l+1}(z) = \frac{z^{l}}{l!}\log
      \left(\frac{1-q^{z}}{1-q}\right)\\
    & &\\  
    & &\qquad  + \sum_{r=0}^{l}\frac{(z-1)^{l-r}}{(l-r)!}L_{r+1}
      +  \sum_{r=1}^{l}\frac{(-1)^{l}z^{l-r}}{(l-r)!}
      \int_{1}^{z}\frac{\xi^{r}}{r!}
    \frac{q^{\xi}\log q}{1-q^{\xi}} d\xi.\\
    & &
  \end{eqnarray*}
we can show that the first term of (\ref{eqn:3*}) vanishes.
Therefore, $\log G_{n}(z+1;q)$ reads as follows:
\begin{prop}
  Suppose $\Re z > -1$ and $m>n$. Then we have 
\begin{eqnarray*}
    & &\log G_{n}(z+1;q)\\
    & &\qquad = \left\{ {z+1 \choose n}+
      \sum_{r=1}^{n} \frac{B_{r}}{r!} \left( - \frac{d}{dz} \right)^{r-1}
      {z \choose n-1} \right\} \log \left(\frac{1-q^{z+1}}{1-q}\right)\\
    & &\\
    & &\qquad + \sum_{r=1}^{n}\left\{ \left( - \frac{d}{dz} \right)^{r-1}
      {z \choose n-1} \right\} \times \int_{1}^{z+1} \frac{\xi^{r}}{r!}
      \frac{q^{\xi} \log q}{1-q^{\xi}} d\xi\\
    & &\\  
    & &\qquad +  \sum_{j=0}^{n-1}G_{n,j}(z)C_{j}(q)
      + \sum_{r=1}^{m} \frac{B_{r}}{r!} F_{n,r-1}(z;q)
      - R_{n,m}(z;q).
\end{eqnarray*}
\end{prop}


\section{Classical limit} 
Now, let us calculate the classical limit. As $q\to1-0$, we see that for 
$\Re z > -1$,
  \begin{eqnarray*}
    & &\mbox{(i)} \quad  \log \left(\frac{1-q^{z+1}}{1-q}\right) \to
      \log \left(z+1\right),\\
    & &\\  
    & &\mbox{(ii)}\quad\int_{1}^{z+1}\frac{\xi^{r}}{r!}
        \frac{q^{\xi}\log q}{1-q^{\xi}}d\xi \to
        -\int_{1}^{\infty}\frac{\xi^{r-1}}{r!}d\xi
        =- \frac{1}{r!r} \left\{\left( z+1 \right)^{r}-1\right\},\\
    & &\\    
    & &\mbox{(iii)}\quad F_{n,r-1}(z;q) \to F_{n,r-1}(z),\\
    & &\\
    & &\mbox{(iv)}\quad R_{n,m}(z;q) \to R_{n,m}(z),\\
    & &\\
    & &\mbox{(v)}\quad C_{j}(q) \to
       C_{j}:= -\sum_{r=1}^{n+1}\frac{B_{r}}{r!}
       \left(\frac{d}{dt}\right)^{r-1}
       \{t^{j} \log t \}\Bigg\vert_{t=1}\\
    & &\qquad \qquad \qquad \qquad
       + \frac{(-1)^{n}}{(n+1)!}\int_{1}^{\infty}
       \overline{B}_{n+1}(t)\left(
       \frac{d}{dt}\right)^{n+1}
       \left\{t^{j}\log t\right\} dt,\\
     \end{eqnarray*}
where
\begin{eqnarray*}
   & & F_{n,r-1}:=\left(\frac{d}{dt}\right)^{r-1}
      \left\{{-t \choose n-1}
        \log\left(\frac{z+t}{z+1}\right)
        \right\}\Bigg\vert_{t=1},\\
  & &\\
  & &  R_{n,m}(z):=\frac{(-1)^{m-1}}{m!}\int_{1}^{\infty}
    \overline{B}_{m}(t) \left(
    \frac{d}{dt}\right)^{m}
    \left\{{-t \choose n-1}
    \log \left(\frac{z+t}{z+1}\right)\right\}dt.\\
  & &
\end{eqnarray*}
Furthermore, from Vitali's convergence theorem, we can show that 
the convergence is uniform on any compact set in $\{\Re z > -1\}.$\par
The constant $C_{j}$ is relevant to the Riemann zeta-function. Indeed,
we obtain 
  \begin{equation}
     C_{j}=\exp(-\zeta'(-j))-\frac{1}{(j+1)^{2}}
  \label{eqn:Cj}\end{equation}   
by applying the Euler-MacLaurin summation formula.\par
Thus we have proved the existence of the classical limit of $G_{n}(z+1;q)$,
and put it to be 
  $$\tilde{G}_{n}(z+1):=\lim_{q \to 1-0}G_{n}(z+1;q).$$
The uniformity of the convergence ensures that $\tilde{G}_{n}(z+1)$
satisfies
  \begin{eqnarray*}
    & & \mbox{1.} \quad \tilde{G}_{n}(z+1)
      =\tilde{G}_{n-1}(z)\tilde{G}_{n}(z),\\
    & & \mbox{2.} \quad \left(\frac{d}{dz}\right)^{n+1}
      \log \tilde{G}_{n}(z+1) \geq 0
      \quad \mbox{for} \quad z \geq 0,\\
    & & \mbox{3.}\quad\tilde{G}_{n}(1)=1,\\
    & & \mbox{4.}\quad \tilde{G}_{0}(z+1)=z+1.
  \end{eqnarray*}
From the uniqueness of the function satisfying these conditions, we have
  $$\tilde{G}_{n}(z+1)=G_{n}(z+1) \quad \mbox{for}\quad \Re z > -1 .$$
Namely, as $q\to1-0$,
  $$G_{n}(z+1;q) \to G_{n}(z+1) \quad in \quad \{\Re z> -1\}.$$
Using the functional equation $G_{n}(z+1)=G_{n-1}(z)G_{n}(z)$, we can show 
the following theorem.

\begin{thm}
As $q \to 1-0$, $G_{n}(z+1;q)$ converges $G_{n}(z+1)$ uniformly 
on any compact set in the domain ${\bf C}\backslash{\bf Z}_{<0}$,
and
  \begin{eqnarray*}
     & &\log G_{n}(z+1)\\
     & & \quad=  \left\{ {z+1 \choose n}+ \sum_{r=1}^{n}
         \frac{B_{r}}{r!}
         \left( - \frac{d}{dz} \right)^{r-1}
         {z \choose n-1} \right\}
          \log(z+1) \\
      & &\\    
      & &\quad - \sum_{r=1}^{n} \left\{
          \left( - \frac{d}{dz} \right)^{r-1}
          {z \choose n-1}\right\}
          \times \frac{1}{r!r}
          \left\{(z+1)^{r}-1\right\}\\
       & &\\   
       & &\quad-  \sum_{j=0}^{n-1}G_{n,j}(z)
          \left\{\zeta'(-j)+\frac{1}{(j+1)^{2}}\right\}
         + \sum_{r=1}^{\infty} \frac{B_{r}}{r!} F_{n,r-1}(z)\\
       & &\\  
       & &\quad - R_{n,m}(z).
    \end{eqnarray*}
\end{thm}

 
\section{The higher Stirling formula}
As $|z|\to\infty$  in $\Delta_{\delta}$, we can see that 
   $$|R_{n,m}(z)| = O(z^{n-m+1}), \quad |F_{m,r-1}(z)|=O(z^{-r+n}).$$
Thus, we have proved the higher Stirling formula:

\begin{thm}
As $|z|\to\infty$ in the sector $\Delta_{\delta}$, we obtain
  \begin{eqnarray*}
    & &\log G_{n}(z+1)\\
    & & \qquad \sim  \left\{ {z+1 \choose n}+ \sum_{r=1}^{n}
      \frac{B_{r}}{r!} \left( - \frac{d}{dz} \right)^{r-1}
      {z \choose n-1} \right\} \log(z+1) \\
    & &\\  
    & & \qquad - \sum_{r=1}^{n}\left\{
      \left( - \frac{d}{dz} \right)^{r-1}
      {z \choose n-1}\right\}\times
      \frac{1}{r!r} \left\{(z+1)^{r}-1\right\}\\
    & &\\
    & & \qquad -  \sum_{j=0}^{n-1}G_{n,j}(z)\left\{
      \zeta'(-j)+\frac{1}{(j+1)^{2}}\right\}
      + \sum_{r=1}^{\infty} \frac{B_{2r}}{(2r)!} F_{n,2r-1}(z).
      \end{eqnarray*}
\end{thm}

\vspace{12pt}

{\bf Examples of the higher Stirling formula.} Let us show some examples.
In the case that $n=1$, this formula coincides with the Stirling formula.
In each case that $n=2$, it coincides with (\ref{eqn:asG}). Namely, we have
  \begin{eqnarray*}
     & &\log G_{2}(z+1)\\
     & & \quad \sim
       \left(\frac{z^{2}}{2}-\frac{1}{12}\right) \log(z+1)
       -\frac{3}{4}z^{2}-\frac{z}{2}+\frac{1}{4} 
       - z \zeta'(0)+ \zeta'(-1) \qquad \qquad \quad\\
     & &\\  
     & & \quad -\frac{1}{12}\frac{1}{z+1}
       + \sum_{r=2}^{\infty}\frac{B_{2r}}{[2r]_{3}}
       \frac{1}{(z+1)^{2r-1}}(z-2r+1).
  \end{eqnarray*}
In the case that $n=3$ and $n=4$, we obtain
  \begin{eqnarray*}
     & &\log G_{3}(z+1)\\
     & & \quad \sim
       \left(\frac{z^{3}}{6}-\frac{z^{2}}{4}+\frac{1}{24}\right)\log(z+1)
       - \frac{11}{36}z^{3}+\frac{5}{24}z^{2}
       + \frac{z}{3} - \frac{13}{72}\\
     & &\\  
     & & \quad - \frac{z^{2}-z}{2} \zeta'(0)
       + \frac{2z-1}{2} \zeta'(-1) -\frac{1}{2}\zeta'(-2)\\
     & &\\
     & &\quad  +\frac{1}{12}\frac{1}{z+1}
       +  \sum_{r=2}^{\infty}
       \left\{z^{2} - (6r-11)z+(4r^2-16r+16)\right\}.\\
& &\\
& &\\
   & &\log G_{4}(z+1)\\
   & &\quad \sim
     \left(\frac{z^4}{24}-\frac{z^3}{6}+\frac{z^2}{6}
     -\frac{19}{720}\right)\log(z+1)\\
   & &\\  
   & &\quad - \frac{4}{72}z^4+\frac{2}{9}z^3
     + \frac{z^2}{8}-\frac{11}{36}z+\frac{31}{144}\\
   & &\\  
   & &\quad - \frac{z^3 - 3z^2 +2z}{6}\zeta'(0)
     + \frac{3z^2-6z+2}{6}\zeta'(-1)
     -\frac{z-1}{2}\zeta'(-2)+\frac{1}{6}\zeta'(-3)\\
   & &\\  
   & &\quad - \frac{1}{12}\frac{1}{z+1}+
     \frac{1}{720}\frac{1}{(z+1)^3}\left(
     6z^2+\frac{13}{2}z+\frac{5}{2}\right)\\
   & &\\
   & &\quad + \sum_{r=3}^{\infty}\frac{B_{2r}}{[2r]_{5}}
     \frac{1}{(z+1)^{2r-1}}\left\{
     z^3-(12r-27)z^2+(20r^2-94r+111)z\right.\\
   & &\\  
   & &\left.\qquad \qquad
     -(8r^3-56r^2+134r-109)
     \right\}.\\
 \end{eqnarray*}


\section{The Weierstrass product representation of the multiple 
gamma function}
In this section, we derive the infinite product representation of the
multiple gamma function more explicitly than (\ref{eqn:vig}).
First, we observe the following proposition.
\begin{prop}
  $$ \exp(\zeta'(-j)) = \exp(P_{j}(1))
     \prod_{k=1}^{\infty}\left\{
      \left(1+\frac{1}{k}\right)^{\frac{B_{r}(k+1)}{j+1}}
      \exp\left(P_{j}(k+1)-P_{j}(k)\right)\right\}$$
where
  \begin{eqnarray*}
    & & P_{j}(x) := \sum_{r=0}^{j+1}\frac{B_{r}}{r!}
      \varphi_{j.r} x^{j-r+1}\\
    & &\\  
    & & \varphi_{j,r} := \left(\frac{d}{dt}\right)^{r}\left\{
      \frac{t^{j+1}}{j+1}\log t - \frac{t^{j+1}}{j+1}
      \right\}\Bigg\vert_{t=1}\\
  \end{eqnarray*} 
and the infinite product converges absolutely.
\end{prop}
   
This proposition is proved by using the Euler-MacLaurin summation formula and
(\ref{eqn:Cj}).\par
Similar calculation shows
\begin{prop}
  \begin{equation}
    \log G_{n}(z+1)=\sum_{j=0}^{n-1}G_{n,j}(z)K_{j}(z)
  \label{eqn;31}\end{equation}
where
  \begin{eqnarray*}
    & & K_{j}(z) := \frac{B_{j+1}(z+1)}{j+1}\log(z+1)
       - \zeta'(-j)+ P_{j}(z+1)\\
    & &\\   
    & & \quad+ \sum_{k=1}^{\infty} \left[
         P_{j}(z+k+1)-P_{j}(z+k) + 
         \frac{B_{j+1}(z+k+1)}{j+1}
         \log \left(
         \frac{z+k+1}{z+k}\right)\right].\\
    & &
  \end{eqnarray*}
Furthermore, the infinite sum of the last term is absolute
convergent.
\end{prop}

By using Proposition 4 and careful consideration on $K_{j}(z)$, we can see
  \begin{eqnarray*}
    & & K_{j}(z) = Q_{j}(z)+ \sum_{r=1}^{j}{j \choose r}z^{j-r}
      \zeta'(-r)-\frac{z^{j+1}}{j+1}\gamma\\
    & &\\  
    & & \qquad + \sum_{k=1}^{\infty}\left\{
      -(z+k)^{k}\log \left(1+\frac{z}{k}\right)
      + \sum_{r=0}^{j}{j \choose r}z^{j-r}
        \sum_{l=1}^{j}\frac{(-1)^{l-1}z^{l}}{l} k^{r-l}
      \right\},
  \end{eqnarray*}
where
  \begin{eqnarray*}
    & &Q_{j}(z) := P_{j}(z+1) -\sum_{r=0}^{j}{j \choose r}z^{r}P_{j-r}(1)\\
    & &\\
    & &\qquad \quad + \frac{1}{j+1}\sum_{r=1}^{j+1}{j+1 \choose r}
    B_{j+1-r}(z)\sum_{l=1}^{r}\frac{(-1)^{l-1}z^{l}}{l}.
  \end{eqnarray*}
The infinite sum in this formula converges absolutely.\par
Finally we can derive the following theorem:

\begin{thm}
For $n \in {\bf N}$, we have
  $$  G_{n}(z+1)
      = \exp \left(F_{n}(z) \right)
         \prod_{k=1}^{\infty}\left\{\left(
           1+\frac{z}{k}\right)^{-{-k \choose n-1}}
           \exp \left( \Phi_{n}(z,k) \right)
         \right\},$$
where
\begin{eqnarray*}
  & & F_{n}(z):= \sum_{j=0}^{n-1}G_{n,j}(z)Q_{j}(z)
    + \sum_{r=0}^{n-2}\left[
      \frac{1}{r!}\left(\frac{\partial}{\partial u}\right)^{r}
      {z-u \choose n-1}\right]_{u=0}^{u=z} \times \zeta'(-r)\\
  & &\quad \qquad - \int_{0}^{z} {z-u \choose n-1} du \times \gamma,\\
  & &\\
  & &\Phi_{n}(z,k) := \frac{1}{(n-1)!} \sum_{\mu=-1}^{n-2}
    \left\{\sum_{r = \mu + 1}^{n-1}
    \frac{{}_{n-1}S_r}{r-\mu}z^{r-\mu}\right\}
    (-1)^{\mu+1}k^{\mu}.\\
\end{eqnarray*}
\end{thm}
 
{\bf Examples of the Weierstrass product representation. }
By making use of Theorem 6, the Weierstrass product representation 
of the multiple gamma function is derived explicitly. We give some
examples. In the case that $n=1$ and $n=2$, the results coincide with
the Weierstrass representation of the gamma function and the $G$-function.
In the case that $n=3$ and $n=4$, we obtain
  \begin{eqnarray*}
    & & G_{3}(z+1)\\
    & & \quad = \exp\left\{-\frac{z^{3}}{4}+\frac{z^{2}}{8}
      +\frac{7}{24}z + \zeta'(-1)-\frac{z(z-1)}{2}\zeta'(0)
      -\left(\frac{z^{3}}{6}-\frac{z^{2}}{4}\gamma\right)\right\}\\
    & &\\  
    & & \quad \qquad \times \prod_{k=1}^{\infty}\left[
      \left(1+\frac{z}{k}\right)^{-\frac{k(k+1)}{2}}
      \exp\left\{\left(\frac{z^{3}}{6}-\frac{z^{2}}{4}\right)\frac{1}{k}
      -\left(\frac{z^{2}}{4}-\frac{z}{2}\right)
      +\frac{z}{2}k\right\}\right],\\
& &\\
& &\\
   & & G_{4}(z+1)\\
   & & \quad = \exp\left\{\frac{61}{144}z^{4}+\frac{13}{18}z^{3}
     + \frac{19}{144}z^{2}-\frac{5}{24}z \right.\\
   & &\\
   & & \quad \qquad \left.
     - \frac{z}{2}\zeta'(-2) + \frac{z^{2}-2z}{3}\zeta'(-1)
     - \frac{z^{3}-3z^{2}+2z}{6}\zeta'(0)
     - \frac{z^{4}-4z^{3}+4z^{2}}{24}\gamma\right\}\\
   & &\\  
   & & \quad \qquad\times \prod_{k=1}^{\infty}\left[
     \left(1+\frac{z}{k}\right)^{\frac{k(k+1)(k+2)}{6}}
     \exp\left\{\left(
     \frac{z^{4}}{24}-\frac{z^{3}}{6}+\frac{z^{2}}{6}\right)
     \frac{1}{k}\right.\right.\\
   & &\\  
   & & \quad \qquad \qquad \qquad \qquad \left.\left.
     - \left(\frac{z^{3}}{18}-\frac{z^{2}}{4}-\frac{z}{3}\right)
     + \left(\frac{z^{2}}{12}-\frac{z}{2}\right)k
     - \frac{z}{6}k^{2}\right\}\right].
  \end{eqnarray*}

\newpage

\noindent{\large\bf Acknowledgement.}\par
The first author deeply thanks to Professor S.Zakrzewski for his
hospitality during the workshop at Banach Center.\par
The first author is partially supported by Grant-in-Aid for Scientific
Research on Priority Area 231 ''Infinite Analysis'' and by Waseda 
University Grant for Special Research Project 95A-257.


\end{document}